\documentclass[12pt]{amsart}

\usepackage{amsmath,amsthm,amssymb,xcolor}
\font\teneufm=eufm10
\font\seveneufm=eufm7
\font\fiveeufm=eufm5
\newfam\frakturfam

\textfont\frakturfam=\teneufm
\scriptfont\frakturfam=\seveneufm
\scriptscriptfont\frakturfam=\fiveeufm

\theoremstyle{definition}
\newtheorem{definition}{Definition}[section]

\theoremstyle{plain}
\newtheorem{lemma}[definition]{Lemma}
\newtheorem{prp}[definition]{Proposition}
\newtheorem{corollary}[definition]{Corollary}
\newtheorem{theorem}[definition]{Theorem}
\newenvironment{Proof}[1][Proof.]{\begin{trivlist}
		\item[\hskip \labelsep {\bfseries #1}]}{\flushright
		$\Box$\end{trivlist}}

\begin{document}

\noindent{\Large
On the generalized Poisson and transposed Poisson  algebras}\footnote{Corresponding author: Farukh Mashurov   (f.mashurov@gmail.com)}

 \bigskip

\begin{center}

 {\bf
Askar Dzhumadil'daev\footnote{Institute of Mathematics and Mathematical Modeling, Almaty, Kazakhstan; \ dzhuma@hotmail.com}, Nurlan Ismailov\footnote{Astana IT University, Astana, Kazakhstan; \ nurlan.ismail@gmail.com}, Farukh Mashurov\footnote{SDU University, Kaskelen, Kazakhstan;}\footnote{Shenzhen International Center for Mathematics, Southern University of Science and Technology, Shenzhen, 518055, China; \ 	f.mashurov@gmail.com
}

}

\end{center} 

 \medskip 
 
\noindent{\bf Abstract}: {\it   We provide the polynomial identities of algebras that are both generalized Poisson algebras and transposed Poisson algebras. We establish defining identities via single operation for generalized Poisson algebras and prove that Itô’s theorem holds for generalized Poisson algebras.}

 \medskip


 \bigskip
 


 \medskip

\section{ Introduction}

 Poisson algebras play an important  role due to their connections with differential geometry and mathematical physics \cite{Arnold1978, Lichnerowicz1977, Weinstein1979}. 
A triple $(A,\cdot,[,])$ is called {\it a Poisson algebra} if $(A,\cdot)$ is a commutative associative algebra, and $(A,[,])$ is a Lie algebra that satisfies {\it the Leibniz rule}:
  \begin{equation}\label{lebniz ident}  [a\cdot b,c]=a\cdot[b,c]+[a,c]\cdot b,\end{equation}
for all elements $a,b,c\in A$. The theory of Poisson algebras has also affected the development of other related algebraic structures, such as noncommutative Poisson algebras \cite{Xu1994}, Jordan dot-bracket algebras \cite{KM1992}, Jacobi algebras, and more generally, generalized Poisson algebras \cite{Agore2015, CK2007, Kirillov1976, Lichnerowicz1977}, and transposed Poisson algebras \cite{BBGW2023}. 

If the Leibniz rule \eqref{lebniz ident}  in the definition of Poisson algebra is replaced by the following identity or also referred to as a compatibility condition:
\begin{equation}\label{lebniz ident with D} 
[a \cdot b, c] = a\cdot [b, c]+[a, c]\cdot b + a\cdot b\cdot D(c),
\end{equation}
then $(A,\cdot,[,],D)$ is called {\it a generalized Poisson algebra},  where $D$ is a derivation of the algebras $(A,\cdot)$ and $(A,[,])$. Note that when $D=0$, the algebra reduces to a Poisson algebra.

The unital generalized Poisson superalgebras appear in the study of simple Jordan superalgebras in the papers \cite{CK2007} and \cite{MarZel2019}. The work of N. Cantarini and V. Kac has significantly contributed to understanding these algebras. They \cite{CK2007} classified all linearly compact simple Jordan superalgebras over an algebraically closed field of characteristic zero. As a consequence, they derived the classification of all linearly compact unital simple generalized Poisson superalgebras. A. Agore and G. Militaru \cite{Agore2015} considered the structure of Jacobi and Poisson algebras, introducing representations of Jacobi algebras, characterizing Frobenius Jacobi algebras, and developing cohomological objects to classify extensions and deformations of these algebraic structures.

Recently, a dual notion of the Poisson algebra, so-called the transposed Poisson algebra, was introduced in \cite{BBGW2023}. This structure arises by exchanging the roles of the two binary operations in the Leibniz rule that defines a classical Poisson algebra.

A triple $(A, \cdot, [,])$ is called {\it a transposed Poisson algebra} if $(A, \cdot)$ is a commutative associative algebra and $(A, [,])$ is a Lie algebra that satisfies the following compatibility identity:
\begin{equation}\label{transpose lebniz ident}  
2a\cdot [b,c] = [a\cdot b,c] + [b,a\cdot c].
\end{equation}
for all elements $a,b,c\in A$.

The theory of transposed Poisson algebras has attracted more attention in recent studies; see \cite{BBGW2023, BFK2024, WeakLeib, Kay2024, Ouar2023, Sar2024} and the references therein. In particular, the authors in \cite{BFK2024} provided an overview of known results in this area and a list of open questions.

In this paper we focus on the characterization and interplay of generalized Poisson algebras and transposed Poisson algebras within the associative commutative algebras endowed with a Lie bracket. We establish that an algebra $(A,\cdot,[,], D)$, with a derivation $D$ can be both a generalized Poisson algebra and a transposed Poisson algebra if and only if it satisfies some certain identities (Theorem \ref{main theorem 1}).

In addition, we will see in Corollaries \ref{generaliation of Bai's Theorem} and \ref{ZinContTrans} that commutative associative and Zinbiel algebras with derivations provides a source of algebras that are both  generalized Poisson and transposed Poisson algebras. We consider Zinbiel algebras with a derivation. In \cite{KMS2024}, it was shown that pre-Novikov algebras coincide with the derived variety of Zinbiel algebras. Furthermore, as proven in \cite{KMS2024}, there are algebras from this derived variety (which matches the class of pre-Novikov algebras) that cannot be embedded into a Zinbiel algebra with a derivation. We show that differential Zinbiel algebras gives  left-symmetric algebras under the multiplication $a\star b=D(a)b-aD(b)$. 
We show that an algebra with a symmetrized Zinbiel product and a skew-symmetrized left-symmetric product gives rise to an algebra that is both a generalized Poisson algebra and a transposed Poisson algebra. In particular, one can have that if $\mathbb{C}[x]$ is the polynomial algebra  over complex numbers $\mathbb{C}$, then  $(\mathbb{C}[X],\circ,[,])$ is both a generalized Poisson algebra and a transposed Poisson algebra with respect to the operations $$[f,g]=D(f)\int_0^x g\, dx -D(g)\int_0^x f\, dx, \quad f\circ g=f\int_0^x g dx +g\int_0^x f dx,$$ where $D=\frac{d}{dx}$ and $f,g\in \mathbb{C}[X]$.

In the study of nonassociative algebras, the concepts of polarization and depolarization play a crucial role in understanding their structure and relationships with algebras. Polarization refers to the technique of representing a one-operation algebra $(A, \ast)$ as an algebra $(A,\circ,\bullet)$ with two operations: a commutative multiplication $(A, \circ)$ and a skew-symmetric multiplication $(A,\bullet)$. For elements $x,y \in A$, the symmetric part is defined as:
$x \circ y = \frac{1}{2}(x \ast y + y \ast x)$
and the skew-symmetric part as:
$x \bullet y = \frac{1}{2}(x \ast y - y \ast x).$
Conversely, depolarization involves starting from an algebra $(A, \circ, \bullet)$ with a commutative product $\circ$ and a skew-symmetric multiplication $\bullet$, from which one can define a single nonassociative multiplication $(A, \ast)$ by:
$x \ast y = x\circ y + x \bullet y .$

This approach  was applied to investigate the structure of nonassociative algebras associated with the Poisson algebras \cite{Markl2006}, \cite{GR2008},  transposed Poisson algebras \cite{WeakLeib}.

Using polarization and depolarization methods, we describe the defining identities with single operations for generalized Poisson algebras (Theorem \ref{PolDepor th1}) and for algebras that are both generalized Poisson and transposed Poisson algebras (Theorem \ref{PolDepor th2}).

In 1955, N. Itô’s \cite{Itô’s1955} proved a fundamental theorem in the structure theory of groups. Itô’s theorem for groups is stated as follows: if  $G$ is a group such that $G=A\cdot B$ for two abelian subgroups $A$ and $B$, then $G$ is metabelian. 

Analogues of Itô’s theorem have also been considered for algebras. For example, for associative algebras \cite{Militaru2017}, Lie algebras \cite{BahturinKegel1995}, \cite{Petravchuk1999}, Leibniz algebras \cite{AgoreMilitaru2015}, Novikov algebras \cite{ShestakovZhang2020}, and Poisson algebras \cite{BCZ23}. We prove that Itô’s theorem holds for generalized Poisson algebras as well (Theorem \ref{Itô’s theorem for GP}).

Throughout the paper, the base field $F$ is algebraically closed of characteristic $0$.

\section{Generalized Poisson and Transposed Poisson algebras}\label{sec: 1}

 In this section, we obtain a necessary and sufficient condition in terms of polynomial identities  for the algebra $(A,  \cdot, [ , ], D)$ to be both a generalized Poisson algebra and a transposed Poisson algebra. Then we provide some examples.

\begin{theorem}\label{main theorem 1}
Let $(A, \cdot)$ be a commutative associative algebra and $(A, [ , ])$ be a Lie algebra. If $D$ is a derivation of $(A,  \cdot, [ , ])$, then $(A,  \cdot, [ , ], D)$ is a generalized Poisson algebra as well as a transposed Poisson algebra if and only if  $(A,  \cdot, [ , ], D)$ satisifies the following identities
  \begin{equation}\label{int iden1}
    a\cdot [b,c]=a \cdot D(b)\cdot c-a\cdot b\cdot D(c),
\end{equation}
\begin{equation}\label{int iden2}
    [a\cdot b,c]=D(a)\cdot b\cdot c+a\cdot D(b)\cdot c-a\cdot b\cdot D(c).
\end{equation}
for all $a,b,c\in A$.
\end{theorem}
\begin{Proof} 
Assume that $(A,  \cdot, [ , ], D)$  is a transposed Poisson and a generalized Poisson algebra with a derivation $D$. To check the identity (\ref{int iden1}), we apply the identity \eqref{lebniz ident with D} to the elements $a\cdot D(b)\cdot c $ and $ a\cdot b\cdot D(c)$, and have  
$$a\cdot[b,c]-a\cdot D(b)\cdot c+a\cdot b\cdot D(c)$$ 
$$=a\cdot [b,c]-[a\cdot c, b]+ a\cdot [c, b] + c\cdot [a, b]+[a\cdot b, c]- a\cdot [b, c] - b\cdot[a, c]$$ 
(by anti-commutativity)
$$a\cdot [c,b]+ c\cdot [a, b]- b\cdot[a, c]+[a\cdot b, c] -[a\cdot c, b]$$ 
(by \eqref{transpose lebniz ident} for $a\cdot [c,b], c\cdot [a, b]$ and $ b\cdot [a, c]$)
$$=\frac{1}{2}[a\cdot c, b]+\frac{1}{2}[c, a\cdot b]+\frac{1}{2}[a\cdot c, b]+\frac{1}{2}[a, c\cdot b]-\frac{1}{2}[a\cdot b, c]-\frac{1}{2}[a, b\cdot c]+[a\cdot b, c] -[a\cdot c, b]=0.$$
Thus, the identity \eqref{int iden1} holds in $A$.
To obtain (\ref{int iden2}) we apply the identity \eqref{lebniz ident with D} to the elements $D(a)\cdot b\cdot c, a\cdot D(b)\cdot c$ and $a\cdot b\cdot D(c)$ and have 
$$[a\cdot b,c]-D(a)\cdot b\cdot c-a\cdot D(b)\cdot c+a\cdot b\cdot D(c)$$
$$=[a\cdot b,c]-[b\cdot c,a]+b\cdot [c,a]+c\cdot [b,a]-[a\cdot c,b]+a\cdot [c,b]+c\cdot [a,b]+[a\cdot b,c]-a\cdot [b,c]-b\cdot [a,c]$$
(by anti-commutativity)
$$=2[a\cdot b,c]-[b\cdot c,a]-[a\cdot c,b]-2a\cdot [b,c]-2b\cdot [a,c]$$
(by \eqref{transpose lebniz ident} for $a\cdot [b,c]$ and $b\cdot [a, c]$)
$$=2[a\cdot b,c]-[b\cdot c,a]-[a\cdot c,b]-[a\cdot b,c]-[b,a\cdot c]-[a\cdot b,c]-[a,b\cdot c]=0.$$ So $(A,  \cdot, [ , ], D)$ satisfies the identity (\ref{int iden2}).

Assume that an algebra $(A,  \cdot, [ , ], D)$ with a derivation $D$ satisfies the identities (\ref{int iden1}) and (\ref{int iden2}). To derive (\ref{lebniz ident with D}) from (\ref{int iden1}) and (\ref{int iden2}), we apply the identity (\ref{int iden1}) for $a\cdot [b,c]$ and $[a,c]\cdot b$, and the identity (\ref{int iden1}) for $[a\cdot b,c]$:
$$[a\cdot b,c]-a\cdot [b,c]-[a,c]\cdot b-a\cdot b\cdot D(c)$$
$$=D(a)\cdot b\cdot c+a\cdot D(b)\cdot c-a\cdot b\cdot D(c)-a\cdot D(b)\cdot c+a\cdot b\cdot D(c)-D(a)\cdot b\cdot c+a\cdot b\cdot D(c)-a\cdot b\cdot D(c)=0.$$

To derive (\ref{transpose lebniz ident}) from (\ref{int iden1}) and (\ref{int iden2}), we apply the identity (\ref{int iden1}) for $a\cdot [b,c]$, and the identity (\ref{int iden2}) for $[a\cdot b,c]$ and $[b,a\cdot c]$:
$$a\cdot [b,c]-\frac{1}{2}[a\cdot b,c]-\frac{1}{2}[b,a\cdot c]$$
$$=a\cdot D(b)\cdot c-a\cdot b\cdot D(c)-\frac{1}{2}D(a)\cdot b\cdot c-\frac{1}{2}a\cdot D(b)\cdot c+\frac{1}{2}a\cdot b\cdot D(c)+\frac{1}{2}D(a)\cdot b\cdot c+\frac{1}{2}a\cdot b\cdot D(c)-\frac{1}{2}a\cdot D(b)\cdot c=0.$$

Hence,  $(A,  \cdot, [ , ], D)$ with the identities (\ref{int iden1}) and (\ref{int iden2}) is a generalized Poisson algebra  as well as a transposed Poisson algebra.
\end{Proof}

A simple construction of an algebra that is both a transposed Poisson algebra and a generalized Poisson algebra is given in the following corollary.
\begin{corollary}\label{generaliation of Bai's Theorem}
Let $(A,\cdot, D)$ be a commutative associative algebra with a derivation $D$. If we define the commutator as $[a,b]=D(a)\cdot b-a\cdot D(b)$ for all $a,b\in L$, then $(A,\cdot, [,], D)$ is a transposed Poisson and generalized Poisson algebra, 
\end{corollary}
This corollary gives an alternative proof of the fact given in \cite[Proposition 2.2]{BBGW2023} that  $(A,\cdot, [,])$ is a transposed Poisson algebra.

Another interesting example can be built by so-called Zinbiel (dual Leibniz) algebras.
An algebra is said to be a {\it Zinbiel} algebra \cite{Loday} if it satisfies the following identity
\begin{equation}\label{Zinbiel identity}
 zin(a,b,c)=(ab+ba)c-a(bc)=0.
\end{equation} 

Zinbiel algebras with a derivation were considered in  \cite{KMS2024}. Let  $A$ be a Zinbiel algebra with a derivation $D$ and define on $A$ two bilinear operations, the left $\prec$ and the right $\succ$ multiplications as follows 
\[
a\prec b = aD(b),\quad a\succ b = D(a)b, \quad a,b\in A.
\]
Then  $(A,\prec, \succ )$ is a pre-Novikov algebra \cite{KMS2024}.
Next proposition shows that every pre-Novikov algebra under the product $a\star b = a\succ b -a\prec b$ is a left-symmetric algebra. 

\begin{prp}\label{left-symmetric}
A Zinbiel algebra $A$ with a derivation $D$ with respect to the product $a\star b= D(a) b-a D(b)$ is a left-symmetric algebra, that is, a Zinbiel algebra $A$ satisfies the identity
    $$(a\star b)\star c -a\star (b\star c)=(b\star a)\star c -b\star (a\star c).$$
\end{prp}
\begin{Proof} Using the the definition of the product $\star$ directly one can have
$$(a\star b)\star c -a\star (b\star c)-(b\star a)\star c +b\star (a\star c)$$
$$=zin(D^2(a),b,c)-zin(a,D^2(b),c)-zin(D(a),b,D(c))+zin(a,D(b),D(c))=0$$
in $A$ for all  $a,b,c\in A$.
\end{Proof}


 \begin{corollary}\label{ZinContTrans}
Let $A$ be a Zinbiel algebra  with a derivation $D$. Then $(A,\circ,[,\,], D)$, where $a\circ b =\frac{1}{2} (a b+b a)$ and $[a, b]=a\star b-b\star a$,  is a generalized Poisson and a transposed Poisson algebra.
\end{corollary}
\begin{Proof}
In \cite{Loday}, it is shown that if $A$ is Zinbiel algebra, then the algebra $(A, \circ)$ under the anti-commutator $a\circ b =\frac{1}{2} (a b+b a)$ is commutative and associative.  Since left-symmetric algebras are Lie admissible, by Proposition \ref{left-symmetric} we obtain that $(A,[,\,]   )$ is a Lie algebra. It is easy to see that $D$ is a derivation of $(A, \circ)$ and $(A, [,])$. Furthermore, we observe that $$[a, b]=a\star b-b\star a=D(a)\circ b -a \circ D(b).$$
Then Corollary \ref{generaliation of Bai's Theorem} implies that $(A,\circ,[,\,], D)$ is a generalized Poisson and a transposed Poisson algebra.
\end{Proof}

The algebra of polynomials over complex numbers $\mathbb{C}[X]$ with respect to the multiplication
$$f\ast g=g\int_0^x f dx$$
becomes a Zinbiel algebra \cite{DzhZinqcom}. Then according to Corollary \ref{ZinContTrans}, one can have an interesting example of  generalized Poisson and a transposed Poisson algebra. Namely, $(\mathbb{C}[X], \circ, [,], D)$ is a generalized Poisson and a transposed Poisson algebra, where
$$[f,g]=D(f)\int_0^x g\, dx -D(g)\int_0^x f\, dx, \quad f\circ g=f\int_0^x g dx +g\int_0^x f dx, \quad D=\frac{d}{dx}.$$

\section{Depolarization of generalized Poisson algebras}

In this section, we describe the generalized Poisson algebras via an identity that involves a single multiplication. 

\begin{theorem}\label{PolDepor th1}  Let $(A,\star, D)$ be an algebra with a derivation $D$ that satisfies the identity
\begin{equation}\label{cont admissible ident}
3(a,b,c)=(a\star c)\star b +(b \star c)\star a - (b \star a)\star c - (c\star a)\star b+4(a\cdot b ) \cdot D(c)+2D(a)\cdot (b\cdot c),
  \end{equation} where $(a,b,c)=(a\star  b)\star  c-a\star  (b\star  c)$ and $a\cdot b=\frac{1}{2}(a\star b+b\star a)$. Then its polarization $(A,\cdot,[,],D)$ is a generalized Poisson algebra.

Conversely,  if $(A,\cdot,[,],D)$ is a generalized Poisson algebra, then its depolarization $(A,\star, D)$ satisfies the identity \eqref{cont admissible ident}.
\end{theorem}

\begin{Proof} 
Let us first define some auxiliary polynomials. 
$$\psi(x,y,z)=-3(x,y,z)+(x\star z)\star y +(y\star  z)\star x - (y\star  x)\star z - (z\star x)\star y+4(x\cdot y ) \cdot D(z)+2D(x)\cdot (y\cdot z),$$
$$Jac(x,y,z)=[[x,y],z]+[[y,z],x]+[[z,x],y],$$
$$Asad(x,y,z)=(x\cdot y)\cdot z-x \cdot (y\cdot z),$$
$$G(x,y,z)=[x\cdot y, z] - x\cdot [y, z]-[x, z]\cdot y - (x\cdot y)\cdot D(z).$$

Suppose that $(A,\star, D)$ satisfies the identity \eqref{cont admissible ident}.  Then $\psi(a,b,c) =0$ holds in $A$ for all $a,b,c\in A$. To see that $(A,\cdot,[,],D)$ is a generalized Poisson algebra, we need to show that $$Jac(a,b,c)=0, \quad Asad(a,b,c)=0, \quad G(a,b,c)=0$$ for all $a,b,c\in A$. 

In fact, the polynomials $Jac(x,y,z), Asad(x,y,z)$ and $G(x,y,z)$ can be written as a linear combination of the polynomial of the form $\psi(x,y,z)$. Specifically, 
$$12Jac(x,y,z)=-\psi(x,y,z)+\psi(x,z,y)+\psi(y,x,z)-\psi(y,z,x)-\psi(z,x,y)+\psi(z,y,x),$$
$$12 Asad(x,y,z)=-\psi(x,y,z)-\psi(x,z,y)+\psi(z,x,y)+\psi(z,y,x),$$
$$12 G(x,y,z)=-\psi(x,y,z)+\psi(x,z,y)-\psi(y,x,z)+\psi(y,z,x)-\psi(z,x,y)-\psi(z,y,x).$$

Therefore, $Jac(a,b,c)=0, Asad(a,b,c)=0, G(a,b,c)=0$ for all $a,b,c\in A$, and consequently, $(A,\cdot,[,],D)$ is a generalized Poisson algebra.

Now suppose that $(A, \cdot, [,], D)$ is a generalized Poisson algebra. Expand  $\psi(x,y,z)$ via $x\star y=x\cdot y+[x,y]$ and have $$\psi(x,y,z)$$$$=-3((x\star y)\star z-x\star(y \star z))+(x\star z)\star y +(y \star z)\star x - (y \star x)\star z - (z\star x)\star y+4(x\cdot y ) \cdot D(z)+2D(x)\cdot (y\cdot z)$$
$$=-3((x\cdot y)\cdot z+[x\cdot y,z]+[x, y]\cdot z+[[x, y], z]-x\cdot(y \cdot z)-x\cdot[y, z]-[x,y\cdot z]-[x,[y, z]])$$
$$+(x\cdot z)\cdot y+[x\cdot z, y] +[x, z]\cdot y + [[x,z], y]+(y\cdot z)\cdot x + [y\cdot z,x]+[y, z]\cdot x+ [[y, z], x]$$
$$-(y\cdot x)\cdot z - [y\cdot x, z] - [y, x]\cdot z- [[y,x],z]-(z\cdot x)\cdot y - [z\cdot x,y] - [z, x]\cdot y- [[z,x],y]$$
$$+4(x\cdot y ) \cdot D(z)+2D(x)\cdot (y\cdot z).$$
Taking into account associativity and the Jacobi identities, one can write 
$$\psi(x,y,z)=-4[x\cdot y, z]-2[x, y]\cdot z+4 [y, z]\cdot x+2[x,y\cdot z]+2[x, z]\cdot y  +4(x\cdot y ) \cdot D(z)+2D(x)\cdot (y\cdot z)$$
$$=-4G(x,y,z)-2G(y,z,x).$$
Then the  identity  $G(a,b,c)=0$ implies that $\psi(a,b,c)=0$ for all $a,b,c\in A$. \end{Proof}

In the case of the trivial derivation $D=0$, note that the identity (\ref{cont admissible ident}) describes Poisson algebra through a single multiplication found in \cite[Proposition 1]{GR2008}.

\section{Generalized Poisson and transposed Poisson admissible algebras}

In this section, we obtain the defining identities of both generalized Poisson and transposed Poisson  algebras in a single multiplication.

Let us first set the following  polynomials:
$$Int_1(x,y,z)=x\cdot [y,z]-(x\cdot z)\cdot D(y)+(x \cdot y) \cdot D(z),$$ 
$$Int_2(x,y,z)= [x  \cdot y, z]- (y\cdot z)\cdot D(x)-(x\cdot z) \cdot D(y) +(x\cdot y)\cdot D(z),$$
$$ f(x,y,z)=x\star  (y\star z)-(z\star y)\star x-2(y\cdot z)\cdot D(x)+4(x\cdot y)\cdot D(z),$$
$$g(x,y,z)=(x\star y)\star z-x\star (y\star z)-\frac{1}{2}(x\star z)\star y+\frac{1}{2}(z\star x )\star y-2(x\cdot y)\cdot D(z),$$ where $x\cdot y= \frac{1}{2}(x\star y + y\star x)$ and 
$[x,y]=\frac{1}{2}(x\star y-y\star x).$

\begin{lemma}\label{lem: Assoc adm} Let $(A, \star, D)$ be an algebra with the identities \begin{equation}\label{contmidcom id1}
    a\star(b\star c)-(c\star b)\star a-2(b\cdot c)\cdot D(a)+4(a\cdot b)\cdot D(c)=0,
\end{equation}
\begin{equation}\label{contmidcom id2} (a\star b) \star c-a\star(b\star c)-\frac{1}{2}(a\star c)\star b+\frac{1}{2}(c\star a )\star b-2(a\cdot b)\cdot D(c)=0.\end{equation}
where $a\cdot b= \frac{1}{2}(a\star b + b\star a).$ 
Then $(A, \cdot)$ is an associative algebra, and $(A, [,])$ is a Lie algebra, where $[a,b]= \frac{1}{2}(a\star b-b\star a)$.
\end{lemma}

\begin{Proof} Suppose that $(A, \star, D)$ satisfies the identities (\ref{contmidcom id1}) and (\ref{contmidcom id2}). Then by $a\cdot b=\frac{1}{2}(a\star b+b\star a)$  we expand
$$12((x\cdot y)\cdot z-x \cdot (y\cdot z))$$
$$=3((x\star y)\star z+(y\star x)\star z+z\star (x\star y)+z\star(y\star x)-x\star (y\star z)-x\star (z\star y)-(y\star z)\star x-(z\star y)\star x)$$
$$=5 f(x,y,z)+f(x,z,y)- f(z,x,y)+3 f(z,y,x)+8 g(x,y,z)+4 g(x,z,y)-4 g(z,x,y).$$
As $f(a,b,c)=0$ and $g(a,b,c)=0$ in $  (A, \star, D),$ for all $a,b,c\in A$, we obtain $$(x\cdot y)\cdot z-x \cdot (y\cdot z)=0,$$ for every $x,y,z\in A.$
For the second part of the lemma we have 
$$4Jac(x, y,z)=-4Asad(z, x, y)-   f(x, y,z)+ 
  f(x,z, y)- 2  g(y, x,z)+2 g(z, x, y).$$

Since $f(a,b,c)=0,$ $g(a,b,c)=0$ and $Asad(a,b,c)=0$ in $(A, \star, D),$ for every $a,b,c\in A$ we obtain the desired result.

 \end{Proof}

\begin{theorem}\label{PolDepor th2}  Let $(A,\star, D)$ be an algebra with a derivation $D$ that satisfies the identities \eqref{contmidcom id1} and \eqref{contmidcom id2}. Then its polarization $(A,\cdot,[,],D)$ is a both generalized Poisson algebra and a transposed Poisson algebra.

Conversely,  if $(A,\cdot,[,],D)$ is a both generalized Poisson algebra and a transposed Poisson algebra, then its depolarization $(A,\star, D)$ satisfies the identity (\ref{contmidcom id1}) and (\ref{contmidcom id2}).
\end{theorem}

\begin{Proof} Suppose that $(A, \cdot, [,], D)$  is a both generalized Poisson algebra and a transposed Poisson algebra. Set $a \star b = a\cdot b + [a,b]$ and have
$$a\star  (b\star c)-(c\star b)\star a-2(b\cdot c)\cdot D(a)+4(a\cdot b)\cdot D(c)$$
$$=a\cdot (b\cdot c)+a\cdot [b, c]+[a, b\cdot c]+[a, [b,  c]]-(c\cdot b)\cdot a-[c,  b]\cdot a-[c\cdot b,a]-[[c, b], a]$$
$$-2(b\cdot c)\cdot D(a)+4(a\cdot b)\cdot D(c)$$
(by \eqref{int iden1} and \eqref{int iden2})
$$=2a\cdot [b,c]-2[c\cdot b,a] -2(b\cdot c)\cdot D(a)+4(a\cdot b)\cdot D(c)=0.$$
So $(A,\star, D)$ satisfies the identity \eqref{contmidcom id1}. 

Now we show that the algebra $(A, \star, D)$ satisfies  the identity \eqref{contmidcom id2}. We have 
$$(a \star b)\star c-a\star (b\star c)-[a,c]\star b-2(a\cdot b)\cdot D(c)$$
$$=(a\cdot b)\cdot c+[a, b]\cdot c+[a\cdot b, c]+[[a, b], c]-a\cdot(b\cdot c)-a\cdot[b, c]-[a, b\cdot c]-[a, [b, c]]$$$$-[a, c]\cdot b-[[a, c],b]-2(a\cdot b)\cdot D(c)$$
(by Lemma \ref{lem: Assoc adm})
$$=[a,b]\cdot c+[a\cdot b,c]-a\cdot[b,c]-[a,b\cdot c]-[a,c]\cdot b-2(a\cdot b)\cdot D(c)$$
$$=Int_1(a, c, b)+Int_1(b, c, a)- Int_1(c, b, a)+
  Int_2(b, a, c)+Int_2(c, b, a).$$

As $Int_1(x, y, z)=0$ and $Int_2(x, y, z)=0$ in $(A, \cdot, [,], D),$  for every $x,y,z\in A,$ we have that the algebra $(A, \star, D)$ satisfies  the identity \eqref{contmidcom id2}. 

Conversely, suppose that $(A,\star, D)$ is an algebra with the identities \eqref{contmidcom id1} and \eqref{contmidcom id2}. Set $a\cdot b= \frac{1}{2}(a\star b + b\star a)$ and 
$[a, b]=\frac{1}{2}(a\star b-b\star a)$. Recall that according to Lemma \ref{lem: Assoc adm} $(A,\cdot, D)$ is an associative algebra and $(A,[,], D)$ is  a Lie algebra. First we will show that $(A,\cdot, [,], D)$ satisfies the identity \eqref{int iden1}. Then we have 
$$a\cdot [b, c]-(a\cdot c)\cdot D(b)+(a \cdot b) \cdot D(c)$$
$$=\frac{1}{4}(a\star(b\star c)-a\star (c\star b)+(b\star c)\star a-(c\star b)\star a)-(a\cdot c)\cdot D(b)+(a \cdot b) \cdot D(c) $$
(by the identity \eqref{contmidcom id1} we have)
$$=\frac{1}{2}(b\cdot c)\cdot D(a)-(a\cdot b)\cdot D(c)-\frac{1}{2}(b\cdot c)\cdot D(a)+(a\cdot c)\cdot D(b)-(a\cdot c)\cdot D(b)+(a \cdot b) \cdot D(c)=0.$$
Therefore, the identity \eqref{int iden1} holds in $(A, \cdot,[,], D)$.

Now we will examine the identity \eqref{int iden2}. We have 
$$[a  \cdot b, c]- (b\cdot c)\cdot D(a)-(a\cdot c) \cdot D(b) +(a\cdot b)\cdot D(c)$$
$$=\frac{1}{4}((a\star b)\star c+(b\star a)\star c-c\star (a\star b)-c\star (b\star a))-(b\cdot c)\cdot D(a)-(a\cdot c) \cdot D(b) +(a\cdot b)\cdot D(c)$$
(by the identity \eqref{contmidcom id1})
$$=\frac{1}{4}(-2(a\cdot b)\cdot D(c)+4(a\cdot c)\cdot D(b)-2(a\cdot b)\cdot D(c)+4(b\cdot c)\cdot D(a))$$$$-(b\cdot c)\cdot D(a)-(a\cdot c) \cdot D(b) +(a\cdot b)\cdot D(c)=0.$$
It implies that the identity \eqref{int iden2} holds in $(A, \cdot, [,], D).$ 
\end{Proof}

Note that a generalized Poisson algebra $(A, \cdot, [,], D)$  with $D(a)=0$, for all $a\in A$, becomes a Poisson algebra. 
An algebra $(A, \cdot, [,])$ is both a Poisson algebra and a transposed Poisson algebra  \cite{BBGW2023} if and only if
\begin{equation}\label{int PtP}
    a\cdot[b,c]= [a\cdot b, c]=0 \text{ for every } a,b,c\in A.
\end{equation}

One can note that in case of $D(a)=0$ in the identities \eqref{contmidcom id1} and \eqref{contmidcom id2}, Theorem \ref{PolDepor th2} provides  identities that define the depolarization of  both Poisson and transposed Poisson algebras, namely, 
$$a\star(b\star c)=(c\star b)\star a,$$
$$(a\star b) \star c-a\star(b\star c)=\frac{1}{2}(a\star c)\star b-\frac{1}{2}(c\star a )\star b.$$
An algebra with the first identity mentioned above is called  {\it reverse associative}  algebra and the free part of reverse associative algebras were studied in the paper \cite{Revass}.

\section{Itô’s theorem for generalized Poisson algebras}

The analogue of Itô’s theorem for Poisson algebras was proved in  \cite{BCZ23}. In this section, we prove that Itô’s theorem also holds for generalized Poisson algebras.

\begin{definition}[\cite{BCZ23}] An algebra $(P,\cdot,[,])$ is called abelian or trivial if both operations are trivial, that is, $x \cdot y =[x, y]=0$ for all $x$ and $y$ in $P.$ \end{definition}

According to the general definition of metabelian groups and Lie algebras the authors \cite{BCZ23} presented a definition of metatrivial Poisson algebra. Following them we rewrite it for generalized Poisson algebras $(P,\cdot,[,]).$

\begin{definition} Let $(P,\cdot,[,])$ be a generalized Poisson algebra.  Then $(P,\cdot,[,])$  is called metatrivial (or strongly metabelian) if $P$ is an extension of an abelian generalized Poisson algebra by another abelian generalized Poisson algebra. \end{definition}

In fact, the authors \cite{BCZ23} established Itô’s theorem for non-commutative Poisson algebras. We will also demonstrate Itô’s theorem for non-commutative generalized Poisson algebras.

\begin{theorem}\label{Itô’s theorem for GP}
\textbf{(Itô’s theorem for non-commutative generalized Poisson algebras)} Let $P$ be a non-commutative generalized Poisson algebra. If $A$ and $B$ be two abelian subalgebras of $P$ such that  $P=A+ B$, then $P$ is metatrivial. \end{theorem}
\begin{Proof}
Let $c_1,\ldots, c_4\in P$ and $c_i=a_i+b_i$ for $1\leq i \leq 4,$ where $a_i \in A$ and $b_i \in B$. For our convenience, we will write $ab$ for the commutative associative product $a\cdot b$. 
Let us show that $c_1c_2c_3c_4=0$ for all $c_1,\ldots, c_4\in P.$ We have   
$$c_1c_2c_3c_4=a_1b_2a_3b_4+b_1a_2b_3a_4.$$ 
Assume that $b_2a_3=a_5+b_5$ and $a_2b_3=a_6+b_6$ for some  $a_5,a_6\in A$ and $b_5,b_6\in B$. Then 
$$c_1c_2c_3c_4=a_1b_2a_3b_4+b_1a_2b_3a_4=a_1a_5b+a_1b_5b_4+b_1a_6a_4+b_1b_6a_4=0.$$
Therefore, for all $c_1,\ldots, c_4\in P,$ we have \begin{equation}\label{itoas}
    c_1c_2c_3c_4=0.
\end{equation}

Next we show that
$c_1[c_2, c_3]c_4= 0.$ We have 
$$c_1[c_2,c_3]c_4= a_1[a_2, b_3]a_4+ a_1[a_2, b_3]b_4+ a_1[b_2, a_3]a_4+ a_1[b_2, a_3]b_4$$
$$+b_1[a_2, b_3]a_4+ b_1[a_2, b_3]b_4+ b_1[b_2, a_3]a_4+ b_1[b_2, a_3]b_4.$$
First we will show that $a_1[a_2, b_3]a_4=0$ and $a_1[a_2, b_3]b_4=0$. 
By  \eqref{lebniz ident with D} one can have
$$a_1[a_2, b_3]a_4=[a_1a_2, b_3]a_4-[a_1, b_3]a_2a_4-a_1a_2D(b_3)a_4= 0.$$
Assume $[a_2, b_3]=a_5+b_5$, for some $a_5\in A$ and $b_5\in B$. Then
$$a_1[a_2, b_3]b_4= a_1a_5b_4 + a_1b_5b_4= 0.$$
Thus one can easily establish that 
\begin{equation}\label{rel1}
    c_1[c_2, c_3]c_4= 0.
\end{equation}

By the identities \eqref{lebniz ident with D},\eqref{itoas} and \eqref{rel1}, we have
$$a_1b_2[b_4, a_3]=a_1[b_2b_4, a_3]-a_1[b_2, a_3]b_4-a_1b_2b_4D(a_3)= a_1[a_3, b_2]b_4= 0,$$
$$[b_4, a_3]b_2a_1=[b_2b_4, a_3]a_1-b_4[b_2, a_3]a_1-b_2b_4D(a_3)a_1= a_1[a_3, b_2]b_4= 0,$$
$$b_1a_2[a_3, b_4]=b_1[a_2a_3,b_4]-b_1[a_2,b_4]a_3-b_1a_2a_3D(b_4)=0,$$
$$[a_3, b_4]a_2b_1=[a_2a_3,b_4]b_1-[a_2,b_4]a_3b_1-a_2a_3D(b_4)b_1=0$$
They imply that
$$c_1c_2[c_3, c_4]=a_1b_2[a_3, b_4]+a_1b_2[b_3, a_4]+b_1a_2[a_3, b_4]+b_1a_2[b_3, a_4]$$
$$=-a_1b_2[b_4,a_3]+a_1b_2[b_3,a_4]+b_1a_2[a_3,b_4]-b_1a_2[a_4,b_3]=0.$$
In a similar way, we can have $[c_1, c_2]c_3c_4= 0$. So by the identity \eqref{lebniz ident with D}, we obtain 
$$[c_1c_2, c_3c_4]=c_1[c_2, c_3c_4]+[c_1, c_3c_4]c_2+c_1c_2D(c_3)c_4+c_1c_2c_3D(c_4)$$
$$=c_1[c_2, c_3]c_4+ c_1c_3[c_2, c_4]-c_1c_3c_4 D(c_2)+c_3[c_1, c_4]c_2+[c_1, c_3]c_4c_2-c_2c_3c_4 D(c_1)= 0.$$

Hence
\begin{equation}\label{sumuprel2}
    c_1c_2[c_3, c_4]=[c_1, c_2]c_3c_4=[c_1c_2, c_3c_4]=0.
\end{equation}

Applying first the identity \eqref{lebniz ident with D} and then by \eqref{itoas}, \eqref{rel1} and \eqref{sumuprel2}, we have
$$[c_1c_2c_3,c_4]=c_1c_2[c_3,c_4]+[c_1c_2,c_4]c_3+c_1c_2c_3D(c_4)$$
$$=c_1c_2[c_3,c_4]+c_1[c_2,c_4]c_3+[c_1,c_4]c_2c_3+c_1c_2D(c_4)c_3+c_1c_2c_3D(c_4)=0.$$

Note that
$$[c_1c_2, [c_3, c_4]]=[a_1b_2, [a_3, b_4]]+[b_1a_2, [a_3, b_4]]+[a_1b_2, [b_3, a_4]] + [b_1a_2, [b_3, a_4]].$$
In order to show that $[c_1c_2, [c_3, c_4]] = 0,$ it suffices to show $[a_1b_2, [a_3, b_4]] = 0.$  
We first show that $a_1b_2D([a_3,b_4])=0.$ By \eqref{sumuprel2} we have
$$a_1b_2D([a_3,b_4])=D(a_1b_2)[a_3,b_4]-D(a_1b_2[a_3,b_4])=D(a_1b_2)[a_3,b_4]$$$$=D(a_1)b_2[a_3,b_4]+a_1D(b_2)[a_3,b_4]=0.$$
The following can be derived in a similar way \begin{equation}\label{rel3}
    b_2a_1D([a_3,b_4])=0.
\end{equation} 
Now assume that $[b_2, a_3]=a_5+ b_5$ and  $[a_1, b_4]=a_6+ b_6$ for some $a_5,a_6\in A$, $b_5,b_6\in B$. Then, by the identity \eqref{lebniz ident with D} and the Jacobi identity, we have
$$[a_1b_2,[a_3,b_4]]=a_1[b_2,[a_3,b_4]]+[a_1,[a_3,b_4]]b_2-a_1b_2D([a_3,b_4])$$
$$=a_1[[b_2, a_3], b_4]+a_1[a_3, [b_2, b_4]] + [[a_1, a_3], b_4]b_2+[a_3, [a_1, b_4]]b_2$$
$$= a_1[a_5, b_4]+[a_3, b_6]b_2$$
$$=[a_1a_5, b_4]-[a_1, b_4]a_5-a_1a_5D(b_4)-[b_6b_2,a_3]-b_6[a_3, b_2]-b_6b_2D(a_3)$$$$=-b_6a_5+ b_6a_5=0.$$

It is clear that to have $[c_1, c_2][c_3, c_4]=0 ,$ it is sufficient to show that  $[a_1, b_2][a_3, b_4]=0$.
Let us again assume that $b_2a_3= a_7+ b_7$ and  $[a_1, b_4]=a_6+ b_6.$
Then by the identity \eqref{lebniz ident with D}, the Jacobi identity  and \eqref{rel3}, we obtain
$$[a_1, b_2][a_3, b_4]=[[a_1, b_2]a_3, b_4]-[[a_1, b_2],b_4]a_3-[a_1, b_2]a_3D(a_4)$$
$$= [[a_1, b_2a_3], b_4]-[b_2[a_1, a_3], b_4] -[b_2a_3D(a_1),b_4] -[[a_1, b_4], b_2]a_3$$
$$= [[a_1, b_7], b_4]-[a_6, b_2]a_3$$
$$= [[a_1, b_4], b_7]-[a_6, b_2a_3]+b_2[a_6, a_3]+b_2a_3D([a_1,b_4])$$
$$=[a_6, b_7]-[a_6, b_7]=0.$$
This completes the proof.

\end{Proof}

\textbf{Acknowledgement:} This research  was funded by the  Science Committee of the Ministry of Science and Higher Education of the Republic of Kazakhstan (Grant No. AP14869221). 

\end{document}